\documentclass[11pt]{article}

\usepackage{amssymb}
\usepackage{latexsym}
\usepackage[leqno]{amsmath}

\newtheorem{lemma}{Lemma}
\newtheorem{proposition}{Propostion}
\newtheorem{theorem}{Theorem}

\newenvironment{proof}[1][Proof]
           {\medbreak\noindent {\em #1. \enspace}}
           {\enspace $\square$ \par \medbreak}
\newenvironment{definition}[1][{\bf \rmfamily Definition}]
           {\medbreak\noindent {\em #1. \enspace}}
           {\par \medbreak}
\newenvironment{remark}[1][{\bf \rmfamily Remark}]
           {\medbreak\noindent {\em #1. \enspace}}
           {\par \medbreak}
\newenvironment{remarks}[1][{\bf \rmfamily Remarks}]
           {\medbreak\noindent {\em #1. \enspace}} {\par \medbreak}
\makeatletter
\@addtoreset{conjecture}{section}
\@addtoreset{equation}{section}
\@addtoreset{lemma}{section}
\@addtoreset{proposition}{section}
\@addtoreset{theorem}{section}
\makeatother

\def\CC{{\mathbb C}}

\def\RR{{\fontsize{12}{14pt}\selectfont\mathbb{R}}}

\newcommand{\gM}{\mathfrak{M}}
\newcommand{\fM}{\mathfrak{M}}

\newcommand{\s}{\sigma}
\renewcommand{\o}{\omega}
\renewcommand{\a}{\alpha}
\renewcommand{\b}{\beta}

\newcommand{\e}{\varepsilon}

\newcommand{\la}{\lambda}

\renewcommand{\t}{\tau}

\newcommand{\La}{\Lambda}

\newcommand{\pt}{\partial}
\newcommand{\ra}{\rightarrow}
\newcommand{\lra}{\longrightarrow}

\renewcommand{\rm}{\textrm}

\newcommand{\cA}{\mathcal{A}}

\newcommand{\cG}{\mathcal{G}}
\newcommand{\cE}{\mathcal{E}}
\newcommand{\cF}{\mathcal{F}}
\newcommand{\cM}{\mathcal{M}}

\newcommand{\cL}{\mathcal{L}}
\newcommand{\cT}{\mathcal{T}}
\newcommand{\cU}{\mathcal{U}}
\newcommand{\cV}{\mathcal{V}}

\def\YM {\mathop{\rm{YM}}\nolimits}
\def\YMH {\mathop{\rm{YMH}}\nolimits}
\def\injrad {\mathop{\rm{injrad}}\nolimits}

\def\div {\mathop{\rm{div}}\nolimits}

\def\Aut {\mathop{\rm{Aut}}\nolimits}
\def\supp {\mathop{\rm{supp}}\nolimits}

\def\Hom {\mathop{\rm{Hom}}\nolimits}
\def\Vol {\mathop{\rm{Vol}}\nolimits}
\def\Ch {\mathop{\rm{Ch}}\nolimits}
\def\rank {\mathop{\rm{rank}}\nolimits}
\def\dist {\mathop{\rm{dist}}\nolimits}

\def\Tr {\mathop{\rm{Tr}}\nolimits}
\def\Pic {\mathop{\rm{Pic}}\nolimits}

\author{By Gang Tian at MIT and Baozhong Yang\footnote{Research partially
supported by NSF Grant DMS-0104163} at Stanford}
\title{Compactification of the moduli spaces of vortices and coupled vortices}
\date{October 19, 2001}
\begin{document}
\maketitle
\begin{abstract}
Vortices and coupled vortices arise from Yang-Mills-Higgs theories and
can be viewed as generalizations or analogues to Yang-Mills
connections and, in particular, Hermitian-Yang-Mills connections. We
proved an analytic compactification of the moduli spaces of vortices
and coupled vortices on hermitian vector bundles over compact K\"ahler
manifolds. In doing so we introduced the concept of ideal coupled
vortices and characterized the singularities of ideal coupled vortices
as well as Hermitian-Yang-Mills connections. 
\end{abstract}

\begin{section}{Introduction}
In the theory of holomorphic vector bundles on compact K\"ahler
manifolds, one remarkable result is the theorem of
Donaldson-Uhlenbeck-Yau on the Hitchin-Kobayashi correspondence
between the stability of a holomorphic vector bundle and the existence
of an irreducible Hermitian-Einstein metric on the bundle. This
enables us to identify the moduli space of stable holomorphic
structures on a complex vector bundle $E$ with the moduli space of
irreducible unitary Hermitian-Yang-Mills connections (HYM connections)
with respect to any fixed Hermitian metric $h$ on $E$. Uhlenbeck's
removable singularity theorem and compactness theorems (\cite{U1},
\cite{U2}) for Yang-Mills connections on four manifolds gave a
compactification of the moduli space of HYM connections on algebraic
surfaces and more generally, of (anti-)self-dual connections any
Riemannian four manifolds. With Uhlenbeck's and Taubes' work as the
analytical basis, Donaldson defined his polynomial invariants for four
manifolds via this compactified moduli space. Through Donaldson's and
other people's work, this idea has been able to provide many deep and
surprising results on the differential topology of four
manifolds. Generalizations of the compactness theorem for Yang-Mills
connections to higher dimensional manifolds have been given in
\cite{Na} and \cite{T}. One of the main results of \cite{T} is a
detailed description of the analytic compactification of the moduli
space of HYM connections on higher dimensional K\"ahler manifolds.

Looking back at their origins in physics, Yang-Mills theories are
important special cases of the more general Yang-Mills-Higgs theories
in physical models which describe interactions of fields and
particles. In a Yang-Mills-Higgs theory people usually study
connections coupled with sections of certain vector bundles, which are
usually called the Higgs objects, and an action which contains, in
addition to the curvature term, the interaction terms between the
connections and the Higgs objects. Yang-Mills-Higgs (YMH) connections
are the critical points for such an action functional. Well-studied
YMH theories include Ginzburg-Landau vortices on $\RR^2$ and monopoles
on $\RR^3$ (See \cite{JT} for an introduction and some fundamental
results. \cite{BBH} contains more recent developments.) For more
examples of Yang-Mills-Higgs theories, we refer the reader to the
survey article \cite{BDGW}. The interactions of geometric analysis,
geometric topology, algebraic geometry and mathmatical physics in the
study of YMH theories make this subject very intriguing.

There is a natural generalization of the vortex equations on $\RR^2$ to
equations on general compact K\"ahler manifolds. Let $(E, h)$ be a smooth
hermitian vector bundle on a compact K\"ahler manifold $(M,
\o)$. Consider an integrable connection $A$ and a section $\phi \in
\Gamma(E)$, the Higgs field. The following equations are called the
{\em vortex equations},
\begin{eqnarray}
\label{eq:1.3} 
&&\bar{\pt}_A \phi = 0, \\
\label{eq:1.4} 
&&\La F_A - \frac i2 \phi \circ \phi^* + \frac i2 \t I_{E} =  0.
\end{eqnarray}
where $\t$ is a real parameter. We shall call a solution $(A, \phi)$
to (\ref{eq:1.3}) and (\ref{eq:1.4}) a $\t$-vortex or simply a
vortex. A $\t$-vortex is the absolute minimum point for the following
YMH functional
$$ \YMH_{\t}(A, \phi) = \int_{M} |F_A|^2 + |d_A \phi|^2 + \frac 14 |\phi
\circ \phi^* - \t I_{E}|^2. $$
A holomorphic pair $(\cE, \phi)$ consists of a holomorphic vector
bundle $\cE$ and a holomorphic section $\phi$ of $\cE$. In \cite{B2}, stable
holomorphic pairs were defined and a Hitchin-Kobayashi type
correspondence between stability of holomorphic pairs and the existence
of irreducible vortices on them were established.

In this paper we are interested in describing the compactification of
the moduli space of vortices and stable pairs. We shall use the
results on the convergence and compactness of pure Yang-Mills
connections in \cite{T} and the removable singularity theorems for
Hermitian-Einstein metrics in \cite{BS}. A key observation we used in
our proofs is a result from O. Garc\'{\i}a-Prada \cite{G1}. This
result identifies the vortex equations on bundles over a K\"ahler
manifold $M$ with the dimensional reduction of the HYM equations under
an $SU(2)$ action on certain associated bundles on the manifold $M
\times \CC P^1$ (see Section 2 for details.) This correspondence
allows us to apply results on HYM connections to study vortices. In
fact, this idea of dimensional reduction has been explored by
E. Witten et al.

A generalization of the vortex equations, the coupled vortex
equations, were  introduced by Garc\'{\i}a-Prada in
\cite{G1}. These equations turn out to be very natural for the setup
of the dimensional reduction mentioned above. Let $(E_1, h_1)$ and
$(E_2, h_2)$ be Hermitian vector bundles on a compact K\"ahler
manifold $(M, \o)$. Consider integrable connections $A_i$ on $(E_i,
h_i)$ ($i=1, 2$) and a section $\phi$ of $\Hom(E_2, E_1)$. The
equations we shall consider are
\begin{eqnarray}
\label{eq:1.5a} 
&&\bar{\pt}_{A_1 \otimes A_2^*} \phi  =  0, \\
\label{eq:1.5b} 
&&\Lambda F_{A_1} - \frac i2 \phi \circ \phi^* + \frac i2 \t I_{E_1} =
0,\\
\label{eq:1.5c} 
&&\Lambda F_{A_2} + \frac i2 \phi^* \circ \phi + \frac i2 \t' I_{E_2} =
0,
\end{eqnarray}
where $A_1 \otimes A_2^*$ is the induced connection on $E_1 \otimes
E_2^*$ and $\t$ and $\t'$ are real parameters. We see that $\t$ and $\t'$ are
related by the Chern-Weil formula and hence there is only one independent
parameter $\t$. The equations (\ref{eq:1.5a}), (\ref{eq:1.5b}) and
(\ref{eq:1.5c}) are called the {\em coupled ($\t$)-vortex equations} and
solutions $(A_1, A_2, \phi)$ of them are called {\em coupled ($\t$)-vortices} on
$(E_1, E_2)$. The coupled vortice equations are also dimensional reductions of
HYM equations on $M \times \CC P^1$ (see Section 2 for details.)

We define $(A_1, A_2, \phi, S, C)$ as an {\em ideal coupled $\t$-vortex}
on hermitian vector bundles $(E_1, h_1)$ and $(E_2, h_2)$ on $M$ if
the singularity set $S$ is a closed subset of $M$ of finite $H^{n-4}$
Hausdorff measure, the triple $(A_1, A_2, \phi)$ is smooth and satisfies the
coupled vortex equations on $M \setminus S$, and $C$ is a holomorphic
chain of codimension 2 on $M$. We can then define the moduli space
$I\cV_{\t}$ of ideal coupled vortices and a natural weak topology on it (see
Section 4.)

Our first main result is
\begin{theorem}[Theorem \ref{th:4.1}]
The moduli space $I\cV_{\t}$ is compact.
\end{theorem}

This compactness theorem easily implies the following compactification
theorems.

\begin{theorem}[Theorem \ref{th:4.2}]
The moduli space $\cV_{\t}$ of coupled $\t$-vortices on hermitian
vector bundles $(E_1, E_2)$ admits a compactification $\bar{\cV}_{\t}$
which is naturally embedded in $I\cV_{\t}$.
\end{theorem}

\begin{theorem}[Theorem \ref{th:4.4}]
The moduli space $V_{\t}$ of $\t$-vortices on a hermitian bundle $E$
over a compact K\"ahler manifold $(M, \o)$ admits a compactification
in the space of ideal coupled $\t$-vortices on $E$ and $L$, where $L$
is the trivial line bundle with the product metric on $M$.
\end{theorem}

We remark here that the compactification and blowup phenomena of
$\t$-vortices appears to be clear only when we embed them into a space
of coupled vortices (see Section 2 and 4 for details.) Using the
Hitchin-Kobayashi type correspondences for vortices and coupled
vortices, these theorems also give the compactification of the moduli space
of corresponding holomorphic objects, i.e., stable pairs and triples
(see Theorem~\ref{th:4.3} and Theorem~\ref{th:4.5}).

Along with the proof of Theorem \ref{th:4.1}, we obtained the
following removable singularity theorems for admissible HYM
connections (see Section 3 for defintion) and ideal coupled vortices,
in which the (non-removable or essential) singularities of HYM
connections and vortices are characterized precisely. Their proofs are
based on Bando and Siu \cite{BS} and the recent work of Tao and Tian
\cite{TT}.

\begin{theorem}
\label{th:1.4}
Let $(A, S)$ be an admissible HYM connection on a hermitian vector
bundle $E$ over $M^m$. Then there exists
$\e_1 = \e_1(m) > 0$ such that if we let
$$ S_0 = \{ x \in S| \lim_{r \ra 0} r^{4-2m} \int_{B_r(x)} |F_{A}|^2
dv \geq \e_1. \}, $$
then $S_0$ is an analytic subvariety of $M$ of codimension $\geq
3$. The holomorphic bundle $E|_{M \setminus S}$ extends to a reflexive
sheaf $\cE$ over $M$ and $S_0$ is the singularity set of $\cE$, i.e.,
the set where $\cE$ fails to be locally free. After a suitable smooth gauge
transformation, the connection $A$ can be
extended to be a smooth connection on $M \setminus S_0$. 
\end{theorem}

\begin{theorem}
\label{th:1.5}
Let $(A_1, A_2, \phi, S, 0)$ be an ideal coupled vortex on hermitian
vector bundles $E_1$ and $E_2$ over $M$. There
exists $\e_2 = \e_2(m, \t) > 0$ such that if we define
$$ S_0 = \{ x \in S | \lim_{r \ra 0} r^{4-2m} \int_{B_r(x)}
e_{\t}(A_1, A_2, \phi) dv \geq \e_2. \}, $$
where $e_{\t}( \cdot, \cdot, \cdot)$ is the YMH action density, then
$S_0$ is an analytic subvariety of $M$ of codimension $\geq 3$.  The
holomorphic bundles $E_j |_{M \setminus S}$ for $j= 1, 2$ extend to
reflexive sheaves $\cE_j$ and $S_0$ is the union of the singularity
sets of $\cE_1$ and $\cE_2$. After a suitable gauge transformation,
the triple $(A_1, A_2, \phi)$ can be extended smoothly over $M
\setminus S_0$.
\end{theorem}

There are some possible directions for further studies. Studies of the
resolution of singularities of stable sheaves might give us a clearer
picture of the boundary points of the compactified moduli spaces.  One
could also define the concepts of stability and semi-stability for
reflexive sheaf pairs or sheaf triples and establish the corresponding
Hitchin-Kobayashi type correspondences for them. The topology of the
compactified moduli space and its relation with the underlying
K\"ahler manifold and with the parameter $\t$ might also be
interesting for studies. We hope to address some of these issues in a
future paper.

In Section 2, we gave an introduction to the vortex equations and
coupled vortex equations on compact K\"ahler manifolds and reviewed
some relevant results, mainly from the references \cite{B1, G1,
G2}. In Section 3, we collected some known results on
Hermitian-Yang-Mills connections. In Section 4, we stated and proved
our compactification theorems. Finally we proved the removable
singularity theorems, Theorem~\ref{th:1.4} and Theorem~\ref{th:1.5}, in
Section 5.
\end{section}

\begin{section}{Vortices, coupled vortices and dimensional reductions}
Assume that $(M, \o)$ is an $m$-dimensional compact K\"ahler
manifold and $(E, h)$ is a hermitian vector bundle on $M$. Let $\cA$
denote the set of all unitary connections on $(E, h)$. We define
$$ \cA^{1,1} = \{ A \in \cA: F_A^{0,2} = 0 \}, $$
the set of integrable connections on $E$. Consider the
following equations for a pair $(A, \phi) \in \cA^{1,1} \times \Gamma(E)$.
\begin{eqnarray}
\label{eq:2.1} 
&&\bar{\pt}_A \phi = 0, \\
\label{eq:2.2} 
&&\La F_A - \frac i2 \phi \otimes \phi^* + \frac i2 \t I_{E} =  0.
\end{eqnarray}
where $\t$ is a real parameter, $\phi^* \in \Gamma(E^*)$ is the dual of
$\phi$ with respect to the the metric $h$ and $\La$ is the contraction
with the K\"ahler form $\o$. In local coordinates, if $\{dz^i,
d\bar{z}^i\}$ is a basis of $T^*M$ and $\o = g_{i\bar{j}}dz^i
d\bar{z}^j$, $(g^{i\bar{j}}) = (g_{i\bar{j}})^{-1}$, then
$$ \La (f_{i\bar{j}}dz^i d\bar{z}^j) = g^{i\bar{j}} f_{i\bar{j}}.$$
The equations (\ref{eq:2.1}) and (\ref{eq:2.2}) are called the {\em vortex
equations }. Since (\ref{eq:2.1}) simply means that $\phi$ is a
holomorphic section of $E$ with respect to the holomorphic structure
defined by $\bar{\pt}_A$, we sometimes call (\ref{eq:2.2}) the vortex
equation. Taking trace of (\ref{eq:2.2}) and integrating, we see that
only for $\t$ such that
$$ \mu (E) = \frac{\deg(E)}{\rank(E)} \leq
\frac{\t \Vol(M)}{4\pi}, $$ the vortex equation can have
solutions. When the equality in above holds, a vortex is given by a
Hermitian-Yang-Mills connection $A$ and $\phi = 0$. 

We define a functional for a pair $(A, \phi) \in \cA \times
\Gamma(E)$,
\begin{equation}
\label{eq:2.3}
\YMH_{\t}(A, \phi) = \int_M |F_A|^2 + |d_A \phi|^2 + \frac 14 | \phi
\otimes \phi^* - \t I_E|^2
\end{equation}
For $(A,
\phi) \in \cA^{1,1} \times \Gamma(E)$, we may compute directly that
(see \cite{B1} for example),
\begin{align*}
&\YMH_{\t}(A, \phi) = \int_M 2 |\bar{\pt}_A \phi|^2 + |\La F_A -\frac
i2 \phi \otimes \phi^* + \frac i2 \t I_{E}|^2 \\
&\qquad + \t \int_M \Tr(i \La
F_A) + \int_M \Tr(F_A^2)\wedge \frac{\o^{m-2}}{(m-2)!}\\
&\quad= \int_M 2 |\bar{\pt}_A \phi|^2 + |\La F_A -\frac
i2 \phi \otimes \phi^* + \frac i2 \t I_{E}|^2 + 2 \pi \t \deg(E) -
8\pi^2 \Ch_2(E).
\end{align*}
Hence the minimum of $\YMH_{\t}$ is the topological quantity $2 \pi \t
\deg(E) - 8\pi^2 \Ch_2(E)$. This minimum is achieved if and only if
$(A, \phi)$ satisfies the vortex equations.

There is another equivalent viewpoint of the vortex equations. We fix a
holomorphic structure on a complex vector bundle $E$ and denote the
obtained holomorphic bundle by $\cE$. Now we allow the hermitian
metric $h$ to vary. For each $h$ there exists a unique unitary
connection which is compatible with the holomorphic structure of
$\cE$. We denote this connection associated to $h$ by $A_h$ and its
curvature by $F_h$. Assume that $\phi$ is a holomorphic section of $\cE$. The
following equation is called the vortex equation for the hermitian
metric $h$,
\begin{equation}
\label{eq:2.4}
\La F_h - \frac i2 \phi \otimes \phi^* + \frac i2 \t I_E = 0.
\end{equation}
It is a standard result that the above two points of views, i.e.,
fixing the metric to consider special unitary connections and fixing
the holomorphic structure to consider special metrics, are
equivalent. We will sketch the idea of this equivalence here (see
Chap. VII, \S 1 of \cite{Ko} for details). We fix a holomorphic bundle
$\cE$ and a holomorphic section $\phi$ here.  Suppose that $\tilde{h}$ is a
hermitian metric satisfying (\ref{eq:2.4}). For any hermitian metric
$h$ on $E$, there exists $g \in \cG^{\CC}$ in the complex linear gauge
group such that $\tilde{h}(s, t) = h(gs, gt), \forall s, t \in
\Gamma(E)$. Suppose that $A$ and $\tilde{A}$ are the connections
associated to $h$ and $\tilde{h}$ respectively. Define $A'$ by
$$ d_{A'} = g \bar{\pt}_{A} g^{-1} + (g^*)^{-1} \pt_{A} g^* = g
d_{\tilde{A}} g^{-1}. $$ Then $A'$ is a unitary connection with
respect to $h$, and $F_{A'} = g F_{\tilde{A}} g^{-1}$, hence $A'$
solves the vortex equation (\ref{eq:2.2}) with respect to the metric
$h$.

Let $\cE$ be a rank $r$ holomorphic vector bundle over $(M, \o)$, and let
$\phi$ be a holomorphic section of $\cE$ and $\t$ be a real
parameter. For background on coherent sheaves and stability, we refer
the reader to Chapter V of \cite{Ko}.
\begin{definition}
The pair $(\cE, \phi)$ is said to be {\em $\t$-stable} if the following
conditions are satisfied,\\
(1) $\mu(\cE') < \hat{\t} = \t \Vol(M)/ (4\pi)$ for every coherent
subsheaf $\cE' \subset \cE$ with $\rank \cE'>0$.\\
(2) $\mu(\cE/\cE') > \hat{\t}$ for every coherent subsheaf $\cE'
\subset \cE$ with $0 < \rank \cE' < r$ and $\phi \in H^0(M, \cE')$.
\end{definition}
We have the following Hitchin-Kobayashi type theorem from Theorem
2.1.6 and Theorem 3.1.1 of \cite{B1}.
\begin{theorem}
\label{th:2.1} 
If there exists a hermitian metric $h$ on $\cE$ which satisfies the
vortex equation (\ref{eq:2.4}), then the bundle splits as $\cE =
\cE_{\phi} \oplus \cE'$ into a direct sum of holomorphic vector
bundles, with $\cE_{\phi}$ containing $\phi$, such that $(\cE_{\phi},
\phi)$ is $\t$-stable and the remaining summands (which together
comprise $\cE'$) are all stable and each of slope $\hat{\t} = \t
\Vol(M) / (4\pi)$.\\
\indent Conversely, if $(\cE, \phi)$ is $\t$-stable, then there exists a 
hermitian metric $h$ on $\cE$ which is a solution of the vortex
equation (\ref{eq:2.4}). 
\end{theorem}
We define the moduli space of $\t$-vortices on a hermitian vector
bundle $(E, h)$  by
$$ V_{\t} = \{ (A, \phi) \hbox{ vortices on } (E, h) \} / \cG, $$ where
$\cG$ is the unitary gauge transformation group on $(E, h)$ and $g(A,
\phi) = (g(A), g \circ \phi)$ for any $g \in \cG$. We define
the moduli space of ($\t$-)stable (holomorphic) pairs by
$$ M_{\t} = \{ [ (\cE, \phi)] : (\cE, \phi)\hbox{ } \t\hbox{-stable pairs with
underlying bundle } E \},$$
where $[ \cdot ]$ means the isomorphism class of holomorphic bundles
with sections. Theorem~\ref{th:2.1} implies that  there is an injection
$M_{\t} \hookrightarrow V_{\t}$. This fails to be a bijection whenever
there is the reducible phenomenon described in the theorem. We let the
exception set of values for $\t$ be
\begin{equation}
\label{eq:2.4b}
\cT = \{ \t | \hat{\t} = \mu(E') \hbox{ for }E' \hbox{ a subbundle
of } E \}.
\end{equation}
If $\t$ is not in $\cT$, then $V_{\t} = M_{\t}$. This identification of the
moduli space of vortices with the moduli space of stable pairs is completely
analogous to the Hitchin-Kobayashi correspondence.

As we have said in the introduction, vortex equations can be
interpretated  as the dimensional reduction of the Hermitian-Yang-Mills
equation on a certain bundle over $M \times S^2$. To better describe and
explore this dimensional reduction, it is natural to introduce a
generalization of the vortex equations.

We assume again that $(M, \o)$ is an $m$-dimensional compact K\"ahler
manifold and $(E_1, h_1)$, $(E_2, h_2)$ are $C^{\infty}$ complex vector
bundles over $M$ with hermitian metrics. Let $\cA_i^{1,1}$ be the
space of integrable unitary connections on $(E_i, h_i)$ for $i =1, 2$.
For a triple $(A_1, A_2, \phi) \in \cA_1^{1,1} \times \cA_2^{1,1} \times
\Gamma(\Hom(E_2, E_1))$, we consider the following equations:
\begin{eqnarray}
\label{eq:2.5a} 
&&\bar{\pt}_{A_1 \otimes A_2^*} \phi  =  0, \\
\label{eq:2.5b} 
&&\Lambda F_{A_1} - \frac i2 \phi \circ \phi^* + \frac i2 \t I_{E_1} =
0,\\
\label{eq:2.5c} 
&&\Lambda F_{A_2} + \frac i2 \phi^* \circ \phi + \frac i2 \t' I_{E_2} =
0.
\end{eqnarray}
where $A_1 \otimes A_2^*$ is the induced connection on $E_1 \otimes
E_2^*$, $\phi^* \in \Gamma (\Hom(E_1, E_2))$ is the adjoint of $\phi$
with respect to $h_1$ and $h_2$, and $\t$ and $\t'$ are real
parameters. By the Chern-Weil
theory, $\t$ and $\t'$ must satisfy the following relation
\begin{equation}
\label{eq:2.5d}
\t \rank E_1  + \t'\rank E_2 = 4\pi \frac{\deg E_1 + \deg
E_2}{\Vol(M)},
\end{equation}
so that there is only one independent parameter $\t$.

We call a triple $(A_1, A_2, \phi)$ a {\em coupled ($\t$)-vortex} if it
satisfies equations (\ref{eq:2.5a}), (\ref{eq:2.5b}) and
(\ref{eq:2.5c}).   Coupled vortices are the absolute minima of the following Yang-Mills-Higgs type functional on $\cA^{1,1}_1 \times \cA^{1,1}_2 \times \Omega^0(E_1 \otimes E_2^*)$,
\begin{align*} 
\YMH_{\t}(A_1, A_2, \phi)& = \int_M  |F_{A_1}|^2 + |F_{A_2}|^2 +
|d_{A_1 \otimes A_2^*} \phi |^2  \\
& \quad + \frac 14 | \phi \circ \phi^* - \t
I_{E_1}|^2 + \frac 14 | \phi^* \circ \phi + \t' I_{E_2}|^2 dv.
\end{align*}
We denote the integrand above by $e_{\t}(A_1, A_2, \phi)$ and call it
the YMH action density for the triple $(A_1, A_2, \phi)$.

Next we describe the dimensional reduction in the setting of triples
(for details and proofs see \S 3 of \cite{G1}). Let $p: M \times S^2 \ra M$,
$q: M \times S^2 \ra M$ be the natural projections. Assume again that
$(E_1, h_1)$ and $(E_2, h_2)$ are hermitian vector bundles on $M$ ,
and $H^{\otimes 2}$ is the degree $2$ line bundle with the standard metric
$h'$ (up to a constant) on $\CC P^1 = S^2$. We shall consider the bundle
\begin{equation*}
F = F_1 \oplus F_2 = p^* E_1 \oplus (p^* E_2 \otimes q^* H^{\otimes
2})
\end{equation*}
on $M \times S^2$ with the induced metric $h = p^* h_1 \oplus p^*h_2
\otimes q^* h_2'$. Consider the left $SU(2)$ action on $M \times S^2$ which
is trivial on $M$ and standard on $S^2$ (i.e., coming from the Hopf
fibration $SU(2) \ra S^2$ and the product structure of $SU(2)$). There
is a natural lift of this action to the bundle $F$, which we shall
describe below.

Note that the total space of $p^*E_1$ is $E_1 \times S^2$, and the
action of $SU(2)$ on the $p^*E_1$ is defined to be trivial on $E_1$
and standard on $S^2$, similarly we define the action of $SU(2)$ on
$p^*E_2$. Recall that $H^{\otimes 2} = SU(2) \times_{\rho_2} \CC $,
where we regard via the Hopf fibration $SU(2)$ as a principle $S^1$ bundle
over $S^2$ and $\rho_2$ is the representation $S^1 \ra S^1 = U(1)$
given by $\rho_2: e^{i\a} \mapsto e^{i2\a}$. This gives a natural
action of $SU(2)$ on $H^{\otimes 2}$ on the left. Since $q^*
H^{\otimes 2} = M \times H^{\otimes 2}$, we require the action of
$SU(2)$ on $q^* H^{\otimes 2}$ to be trivial on $M$ and as above on
$H^{\otimes 2}$. 

Any unitary connection on $(F, h)$ is then of the form
\begin{equation}
\label{eq:2.6}       d_A =  \begin{pmatrix} d_{\tilde{A_1}} & \b\\
                                            -\b^* & d_{\tilde{A_2}}
                            \end{pmatrix} ,
\end{equation}
where $\tilde{A_1}$, $\tilde{A_2}$ are connections on $(F_1,
\tilde{h_1})$ and $(F_2, \tilde{h_2})$ respectively, and $\b \in
\Omega^1(M \times S^2, \Hom(F_2, F_1))$. It can be shown (as in
Prop. 3.5 of \cite{G1}) that for any $SU(2)$-invariant connection $A$ on
$(F, h)$, we have
\begin{equation}
\label{eq:2.6b}
\tilde{A_1} = p^*A_1, \tilde{A_2} = p^*A_2 \otimes q^*A_2',
\end{equation}
where $A_1$ and $A_2$ are connections on $(E_1,
h_1)$ and $(E_2, h_2)$ respectively, and $A_2'$ is the unique
$SU(2)$-invariant connection on $(H^{\otimes2}, h_2')$. We also have
\begin{equation}
\label{eq:2.7}  
\b = p^*\phi \otimes q^* \a,
\end{equation}
where $\phi \in \Omega^0(M, E_1 \otimes
E_2^*)$, and $\a$ is the unique $SU(2)$-invariant element of
$\Omega^1(S^2, H^{\otimes {-2}})$, up to a constant factor. The proofs
of these facts exploit the $SU(2)$-invariance of the objects and
study the restriction to fibers to show that certain components have
to vanish (for details see the proof of Prop. 3.5 in \cite{G1}).

Let $\cA_1$ and $\cA_2$ be the spaces of unitary connections on $(E_1,
h_1)$ and $(E_2, h_2)$ respectively and $\cA^{SU(2)}$ be the space of
$SU(2)$-invariant unitary connections on $(F,h)$. Then after fixing
the choice of the one form $\a$, we have a one-to-one correspondence
between
\begin{equation}
\label{eq:2.8} 
\cA^{SU(2)} \leftrightarrow \cA_1 \times \cA_2 \times
\Gamma(\Hom(E_2, E_1))
\end{equation}
given by $A \leftrightarrow (A_1, A_2, \phi)$ as in (\ref{eq:2.6}),
(\ref{eq:2.6b} and (\ref{eq:2.7}). If we restrict (\ref{eq:2.8}) to
integrable connections, there is a one-to-one correspondence
\begin{equation}
\label{eq:2.8a}
\cA^{1,1, SU(2)} \leftrightarrow \cA_1^{1,1} \times \cA_2^{1,1} \times
\Gamma(\Hom(E_2, E_1)).
\end{equation}
We denote by $\cG^{SU(2)}$ the $SU(2)$-invariant gauge transformation
group on $F$ and $\cG_i$ the unitary gauge transformation groups on
$E_i$ ($i=1, 2$). Using a similar argument to the above we can 
write every $g \in \cG^{SU(2)}$ as
$$ g = \begin{pmatrix}
                g_1 & 0 \\
                0   & g_2
       \end{pmatrix},
$$
for $g_1 \in \cG_1$ and $g_2 \in \cG_2$. Hence there is
a one-to-one correspondence of the gauge groups,
\begin{equation}
\label{eq:2.8b}
\cG^{SU(2)} \leftrightarrow \cG_1 \times \cG_2.
\end{equation}
It is then clear that there is a one-to-one correspondence between the
configuration spaces,
$$
\cA^{SU(2)} / \cG^{SU(2)} \stackrel{1 - 1}{\leftrightarrow} \cA_1
\times \cA_2 \times \Gamma(\Hom(E_2, E_1))/ \cG_1 \times \cG_2,
$$
where $(g_1, g_2)(A_1, A_2, \phi) = (g_1(A_1), g_2(A_2), g_1 \circ
\phi \circ g_2^{-1})$ is the gauge action of $\cG_1 \times \cG_2$ on
triples.

From now on, we shall fix a K\"ahler metric
$$\o_{\s} = \o \oplus \s
\o_{\CC P^1}$$
on $M \times S^2$, where $\o_{\CC P^1}$ is the standard
Fubini-Study metric on $\CC P^1$ so that $\int_{\CC P^1} \o_{\CC P^1} =
1$, and $\s$ is given by the following formula,
\begin{equation}
\label{eq:2.9} 
\s = \frac{2 r_2\Vol(M)}{(r_1 + r_2) \hat{\t} - \deg E_1 - \deg
E_2},
\end{equation}
where $r_1 = \rank E_1$ and $r_2 = \rank E_2$ and $\hat{\t} = \t
\Vol(M)/(4 \pi)$. In what follows, we shall fix the choice of $\a$ such that
\begin{equation}
\label{eq:2.9a}
\a \wedge \a^* = \frac i2 \s \o_{\CC P^1}.
\end{equation}
Define the Yang-Mills functional for a connection $A$ on the hermitian
vector bundle $F$ with respect to
$\o_{\s}$ by
$$ \YM_{\s}(A) = \int_{M \times S^2} |F_A|_{\s}^2
\frac{\o_{\s}^{m+1}}{(m+1)!}, $$
where $| \cdot |_{\s}$ is the norm induced by the hermitian metric on
$F$ and the K\"ahler metric $\o_{\s}$ on $M \times S^2$.

If $A \in \cA^{SU(2)}$ corresponds  to $(A_1,
A_2, \phi) \in \cA_1 \times \cA_2 \times \Gamma(\Hom(E_2, E_1))$
in (\ref{eq:2.8}),
then we have the following lemma relating the energy densities of
$\YMH_{\t}$ and $\YM_{\s}$,
\begin{lemma}
\label{lem:2.1}
\begin{equation}
\label{eq:2.9b}
|F_A|_{\s}^2 = e_{\t}(A_1, A_2, \phi) + c(\t),
\end{equation}
where $c(\t)$ is a constant depending only on $\t$, $\Vol(M)$ and the
degrees of $E_1$ and $E_2$.
\end{lemma}
\begin{proof}
Via the expression (\ref{eq:2.6}), we compute that
$$
F_A = \begin{pmatrix}
               p^*F_{A_1} - \b \wedge \b^* &  d_{A_1 \otimes A_2^*}
               \b\\
                        &  \\
                -(d_{A_1 \otimes A_2^*} \b)^* &  p^*F_{A_2} - 4 \pi i
               \o_{\CC P^1} + \b^* \wedge \b
        \end{pmatrix}.
$$
We also note that with the convention of (\ref{eq:2.9a}),
\begin{align*}
\b \wedge \b^* &= \frac i2 p^*(\phi \circ \phi^*) \o_{\CC P^1},\\
\b^* \wedge \b &= - \frac i2 p^*(\phi^* \circ \phi) \o_{\CC P^1}.
\end{align*}
These together give us (for simplicity we omit the $p^*$'s in the following
equation),
\begin{align*}
|F_A|_{\s}^2 & = |F_{A_1}|^2 + \frac 14 | \phi \circ \phi^*|^2 +
 |F_{A_2}|^2 + \frac 14 | \phi^* \circ \phi - \frac{8\pi}{\s} I_{E_2}
|^2\\
& \quad + 2|\a|_{\s}^2|d_{A_1 \otimes A_2^*} \phi|^2 + 2 | d_{A_2'} \a|_{\s}^2
 |\phi|^2 \\
& = e_{\t}(A_1, A_2, \phi) + ( \frac{16\pi^2}{\s^2} r_2^2 - \frac 14
(\t^2 r_1^2 + {\t'}^2 r_2^2)).
\end{align*}
In the above we used the fact that $d_{A_2'} \a = 0$, $|\a|_{\s}^2 = |
\a \wedge \a^*| = 1/2$, and the identity $\frac {4
\pi}{\s} = \frac{ \t - \t'}{2}$.
\end{proof}

As a corollary of Lemma~\ref{lem:2.1}, we have the following identity,
\begin{equation}
\label{eq:2.10}
\YM_{\s}(A) = \s \YMH_{\t}(A_1, A_2, \phi) + C(\t),
\end{equation}
where $C(\t)$ is a geometrical constant depending on $\t$. 

By direct computation it can be shown (see \S 3 of \cite{G1}) that
under the correspondence (\ref{eq:2.8}), a connection $A$ on $M
\times S^2$ is a Hermitian-Yang-Mills connection with respect to
$\o_{\s}$ if and only if $(A_1, A_2, \phi)$ is a coupled $\t$-vortex,
where $\s$ is related to $\t$ by (\ref{eq:2.9}).

Assume that $\rank E_2 = 1$, we can define a concept of stability for
holomorphic triples $(\cE_1, \cE_2, \phi)$, where $\cE_1$ and $\cE_2$
are holomorphic bundles with underlying topological bundle $E_1$ and
$E_2$ respectively and $\phi$ is a holomorphic section of $\Hom (\cE_2,
\cE_1)$.
\begin{definition}
Assume that $\rank E_2 =1$. We define a triple $(\cE_1, \cE_2, \phi)$
to be {\em $\t$-stable} if  the pair $(\cE_1 \otimes
\cE_2^*, \phi)$ is $\t$-stable.
\end{definition}
We have the following theorem (Theorem 4.33 in \cite{G1}).
\begin{theorem}
\label{th:2.2}
Let $(A_1, A_2, \phi)$ be a coupled $\t$-vortex and let $\cE_1 = (E_1,
\bar{\pt}_{A_1})$ and $\cE_2 = (E_2, \bar{\pt}_{A_2})$. Then 
$\cE_1$ decomposes as a direct sum $\cE_1 = \cE' \oplus \cE''$, such
that $\phi \in \Hom (\cE_2, \cE')$, $(\cE', \cE_2, \phi)$ is
$\t$-stable and $\cE''$ is a direct sum of stable bundles of the slope
$\hat{\t} = \t \Vol(M)/(4 \pi)$.

Conversely, if $(\cE_1, \cE_2, \phi)$ is $\t$-stable, then for any
hermitian metrics $h_1, h_2$ on $E_1, E_2$ respectively, there exists
a solution $(A_1, A_2, \phi)$ to the coupled $\t$-vortex equations
such that $\cE_1 = (E_1, \bar{\pt}_{A_1})$ and $\cE_2 = (E_2,
\bar{\pt}_{A_2})$ and the solution is unique up to unitary gauge
transformations.
\end{theorem}

The following discussion is similar to the case of vortices and stable
pairs. We define the moduli space of coupled $\t$-vortices on hermitian
bundles $E_1$ and $E_2$ by
$$ \cV_{\t} = \{ (A_1, A_2, \phi) \hbox{ coupled vortices on } (E_1,
E_2) \} / \cG_1 \times \cG_2, $$
where $\cG_1$ and $\cG_2$ are the unitary gauge transformation groups
on $E_1$ and $E_2$. If $\rank E_2 = 1$, we can define the moduli space
of stable triples by
$$ \gM_{\t} = \{ [(\cE_1, \cE_2, \phi)] :(\cE_1, \cE_2, \phi)\hbox{ }
\t\hbox{-stable triple with underlying bundles } E_1 \hbox{ and } E_2
\},$$
where $[ \cdot ]$ means the holomorphic isomorphism class.
Theorem~\ref{th:2.2} implies that when $\rank E_2 =1$, there is an
injection $\fM_{\t} \hookrightarrow \cV_{\t}$. This fails to be a
bijection when there is the reducible phenomenon described in the
theorem. We let
\begin{equation}
\label{eq:2.11}
\cT = \{ \t | \hat{\t} = \mu(E') \hbox{ for }E' \hbox{ a subbundle
of } E_1 \}.
\end{equation}
If $\t$ is not in $\cT$, then $\cV_{\t} = \fM_{\t}$.

Next we consider the relation between stable pairs and triples. Let
$L$ be the trivial smooth line bundle over $M$ and let $\fM_{\t}$ be
the moduli space of $\t$-stable triples with the underlying smooth
bundles being $E$ and $L$. Let $\cE$ and $\cL$ be holomorphic
structures on $E$ and $L$. Since by definition, $(\cE, \cL, \phi)$ is
$\t$-stable if and only if $(\cE \otimes \cL^*, \phi)$ is stable, we
have a map $\fM_{\t} \ra M_{\t}$ given by
$$ [(\cE, \cL, \phi)] \mapsto [(\cE \otimes \cL^*, \phi)]. $$
The group of holomorphic bundles supported by $L$, $\Pic^0(M)$, acts
on $\cM_{\t}$ by
$$ (\cE, \cL, \phi) \mapsto (\cE\otimes \cU, \cL \otimes \cU, \phi),
\quad \hbox{ for } \cU \in \Pic^0(M). $$ It is then clear that the
above map gives an identification
\begin{equation}
\label{eq:2.11b} 
M_{\t} \cong \fM_{\t}/\Pic^0(M)
\end{equation}

Finally we consider the relation between the vortex equations and the
coupled vortex equations. We let $(E_1, h_1) = (E, h)$ of rank $r$ and
$(E_2, h_2) = (L, h_0)$, where $L$ is the trivial bundle and $h_0$ is
the product metric. Then a coupled $\t$-vortex $(A, A', \phi)$ on $(E,
L)$ satisfies
\begin{eqnarray}
\label{eq:2.12a} 
&&\bar{\pt}_{A} \phi =  0, \\
\label{eq:2.12b} 
&&\Lambda F_{A} - \frac i2 \phi \circ \phi^* + \frac i2 \t I_{E} =
0,\\
\label{eq:2.12c} 
&&\Lambda F_{A'} + \frac i2 \phi^* \circ \phi + \frac i2 \t' =
0,
\end{eqnarray}
where $\t r + \t' = 4 \pi \deg E / \Vol(M)$. We see that
(\ref{eq:2.12a}) and (\ref{eq:2.12b}) are actually the vortex
equations on $(E, h)$. Hence $(A, \phi)$ is a vortex on $(E, h)$. On
the other hand, if $(A, \phi)$ is a solution to the vortex equations
(\ref{eq:2.12a}) and (\ref{eq:2.12b}) on $(E, h)$, then it is easy to
check that the linear equation (\ref{eq:2.12c}) admits a solution $A'$
(see for example \cite{G2} p.541 - 542). The solution $A'$ is unique if we
fix the holomorphic structures on  $L$ and require  $A'$
to be compatible with the given holomorphic structure. Therefore,  if
we fix the trivial holomorphic structure on $L$, there is an embedding
\begin{equation}
\label{eq:2.13a}
V_{\t} \hookrightarrow \cV_{\t},
\end{equation}
For each holomorphic structure on $L$, there is a unique solution
$A'$. Hence if we take all the holomorphic structures on $L$ into
consideration, then there is an identification
\begin{equation}
\label{eq:2.13} 
V_{\t} \cong \cV_{\t}/\Pic^0(M).
\end{equation}

\end{section}

\begin{section}{Hermitian-Yang-Mills connections}
In this section we collected some results on the blow-up phenomena of
Yang-Mills connections and in particular, some results on
Hermitian-Yang-Mills connections. The main references for this section
are \cite{T} and \cite{BS}.

Assume that $M$ is an $n$-dimensional Riemannian manifold and $E$ is a
 smooth vector bundle on $M$ with a compact structure group $G$. Let
 $\cA$ be the set of all $G$-connections on $E$.\\[1ex]
The Yang-Mills functional $\YM: \cA \ra \RR$ is defined by
$$ \YM(A) = \int_M |F_A|^2 dv. $$ A smooth connection $A$ is called a
{\em Yang-Mills connection} if and only if $A$ is a critical point of
the functional $\YM$, or equivalently, $A$ satisfies the Yang-Mills
equation
$$ d_A^* F_A = 0. $$
\noindent {\it Definition.} 1) $(A,S)$ is called an {\em admissible
 connection} if $S \subset M$, $H^{n-4}(S)=0$, $A$ is a
smooth connection on $ M\setminus S$ and $\YM(A) < \infty$. $S$ is
called the singular set of $A$.\\[1ex]
2) $(A,S)$ is called an {\em admissible
Yang-Mills connection} if $(A, S)$ is an admissible connection and
$d_A^*F_A = 0$ on $M\setminus S$.\\[1ex]
Sometimes, when the singular set $S$ is understood, we simply say $A$
is an admissible connection.

By using variations generated by a vector field on $M$, we can derive the
following {\em first variation formula} for {\em smooth} Yang-Mills
connections (first derived by Price \cite{P}),
\begin{equation}
\label{eq:var}
\int_M |F_A|^2 \div X - 4 \sum_{1\leq i < j \leq n}\langle F_A
( \nabla_{e_i}X, e_j), \,F_A(e_i, e_j) \rangle dv = 0.
\end{equation}
This formula is true for any compactly supported $C^1$ vector field
$X$ on $M$. This motivates the following definition. \\[1ex]
\noindent {\em Definition. } An admissible connection $(A, S)$ is said to be
{\em stationary} if the first variation formula is true for $A$ with
any compactly supported $C^1$ vector field $X$ on $M$.\\[1ex]
\indent Assume that $A$ is a stationary connection. By using a cutoff
of the radial vector field $X = \sum_{i=1}^{n} x_i \frac{\pt}{\pt x_i}
$ in the first variation formula, we can obtain the following
important monotonicity formula first shown by Price \cite{P} (also see
Tian \cite{T} for a more general version). Let $\injrad(x)$ denotes
the injective radius of $x \in M$.
\begin{proposition}
\label{prop:3.1a} 
For any $x \in M$, there exist positive constants $a = a(x)$ and $r_x
< \injrad(x)$ which only depend on the supremum bound of the curvature of
$M$, such that if $0 < \sigma < \rho \leq r_x$, then
\begin{align}
\label{eq:m}
\rho ^{4-n}e^{a \rho^2}& \int_{B_\rho (x)}\left| F_A\right| ^2dv-\sigma
^{4-n}e^{a \sigma^2}\int_{B_\sigma (x)}\left| F_A\right| ^2dv\\
& \geq 4\int_{B_\rho (x)\setminus B_\sigma \left( x\right)
}r^{4-n}\left| \frac \partial {\partial r} \rfloor F_A \right|^2 dv.    
\nonumber
\end{align}
If $M$ is flat, then we may take $a = 0$. 
\end{proposition}

We have the following a priori pointwise estimate for smooth
Yang-Mills connections  by Uhlenbeck and also by Nakajima
\cite{Na}.  It can be proven by using the Bochner-Weitzenb\"ock
formula and a method similar to Schoen's method (IX.4.2 in \cite{SY})
in proving the a priori pointwise estimate for stationary harmonic
maps.
\begin{proposition}
\label{prop:3.1} 
Assume that $A$ is a smooth Yang-Mills connection. There exist $\e_0 = \e_0(n)
> 0$ and $C = C(M, n)> 0$ such that for any $x\in M$, $\rho < r_x$, if
$\rho^{4-n} \int_{B_\rho(x)} |F_A|^2 \, dv \leq \e_0$, then
\begin{equation}
\sup_{B_{\frac 12 \rho}(x)} \rho^2|F_A|(x) \leq C \left( \rho^{4-n} \int_{B_\rho(x)} |F_A|^2 \, dv
\right)^{\frac12}
\end{equation}
\end{proposition}

Using the a priori pointwise estimate and monotonicity formula,
compactness theorems about Yang-Mills connections were first proven in
Uhlenbeck \cite{U2} and Nakajima \cite{Na}. We quote the following
compactness theorem from \cite{T}.
\begin{proposition}
\label{prop:3.2}
Assume that $\{ (A_i,\, S_i)\}$ is a sequence of stationary admissible
Yang-Mills connections on $E$, with $YM(A_i)\leq \Lambda $, where
$\Lambda $ is a constant. Let $S_{cls} = \overline{ \limsup_{i \ra
\infty} S_i}$. Assume also that $H^{n-4}(S_{cls}) = 0$. Then there
exist a subsequence $\{A_{i}\}$, a closed subset $S_b$ of $M$ with
$H^{n-4}\left( S_b \cap K \right) <\infty $ for any compact subset $K
\subset M$,  a nonnegative $H^{n-4}$-integrable
function $\Theta$ on $S_b$, gauge transformations $\sigma_i \in
\Gamma(\Aut P)$ and a smooth Yang-Mills connection on $M \setminus
S_b$, such that the following holds:

(1) On any compact set $K\subset M\backslash ( S_b\cup S_{cls})$, $\sigma_i(A_{i})$ converges to $A$ in $C^\infty $ topology.

(2) $|F_{A_{i_j}}|^2 dv \lra |F_A|^2 dv + \Theta H^{n-4} \lfloor S_b$
    weakly as measures on U.
\end{proposition}
\begin{remarks}
1) We shall make a little remark about the condition $H^{n-4}(S_{cls})
= 0$ here. This condition will be trivially satisfied if all $A_i$ are
smooth YM connections, i.e., all $S_i$'s are empty. This condition is
necessary for the application of the {\em a priori} estimates to
extract a convergent subsequence. We conjecture that the theorem is
true without this condition, but the proof will need a very general
regularity theorem for stationary admissible YM connections, which
has not been proved yet.\\
2) The density function $\Theta$ is defined
by
$$ \Theta (x) = \lim_{\rho \ra 0} \liminf_{j \ra \infty} \rho^{4-n}
\int_{B_{\rho}(x)} |F_{A_{i_j}}|^2 dv. $$
3) The closed set $S_b$ is given by
\begin{equation}
\label{eq:3.1}
S_b  = \{ x \in M : \Theta(x) \geq \e_0 \}.
\end{equation}
where $\e_0$ is as in Prop.~\ref{prop:3.1}.
\end{remarks}
\indent Define
$$ T = \hbox{ the closure of } \{ x \in S_b | \Theta(x)>0 , \lim _{r
\ra 0+} r^{4-n} \int_{B_r(x)} |F_A|^2 dv = 0 \}. $$
It is easy to show that the measue $ \Theta H^{n-4}\lfloor S_b$ is
equal to $\Theta H^{n-4} \lfloor T$.  $(T, \Theta)$ is called the {\em
blow-up locus} of the sequence $\{A_{i_j}\}$. The set $T$ may be shown
to be rectifiable (see \S 3.3 of \cite{T}), and hence $(T, \Theta)$
defines a rectifiable varifold.  We know very little about the blow-up
set of sequences of general Yang-Mills connections without further
restrictions.

The most important examples of Yang-Mills connections include
self-dual and anti-self-dual connections on four manifolds, and
Hermitian-Yang-Mills connections (also called Hermitian-Einstein
connections in literature) on Hermitian vector bundles over K\"ahler
manifold. Assume that $(M, \o)$ is an $m$-dimensional K\"ahler manifold and
$(E, h)$ is a Hermitian vector bundle over $M$. A unitary connection $A$
on $(E, h)$ is a {\em Hermitian-Yang-Mills connection} (HYM
connections) if it is integrable and
$$\La F_A = \la I_E,$$ where $\la$ is a constant. HYM connections are
the absolute minima of the Yang-Mills functional.

One of the main results in \cite{T} (Theorem 4.3.3) is the following
characterization of the blow-up locus of a sequence of HYM
connections (in \cite{T}, a more general class of connections,
$\Omega$-anti-self-dual connections are treated).
\begin{theorem}
\label{th:3.1}
Let $(M, \o)$ be an $m$-dimensional compact K\"ahler manifold and $(E,
h)$ a hermitian vector bundle over $M$. Let $\{A_i\}$ be a sequence of
Hermitian-Yang-Mills connections on $E$. Then by passing to a
subsequence, $A_i$ converges to an admissible Hermitian-Yang-Mills
connection $A$ (in the sense in Prop.~\ref{prop:3.2}) with the blow-up
locus equivalent as a $(2m-4)$ varifold to $(S, \Theta)$, such that $S
= \cup_{\a} S_{\a}$ is a countable union of $(m-2)$ dimensional
holomorphic subvarities $S_{\a}$'s, and $\frac{1}{8\pi^2}
\Theta|_{S_{\a}}$ are positive integers for every $\a$. There is the
following convergence of measures,
$$ |F_{A_i}|^2 dv \ra |F_A|^2 dv + \Theta H^{2m-4}\lfloor S. $$
\end{theorem}
The following removable singularity/extension theorem is proved in
\cite{BS} and the proof is based on results in \cite{B, S}.
\begin{theorem}
\label{th:3.2}
Assume that $(E, h)$ is a holomorphic vector bundle defined outside a
closed subset $S$ with $H^{2m-4}(S) < \infty$, on an $m$-dimensional
K\"ahler manifold $(M, \o)$. Assume that the curvature $F_h$ isn
locally square integrable on $M$, then $E$ extends uniquely to a
reflexive sheaf $\cE$ over $M$. If furthermore, $h$ satisfies the
Hermitian-Einstein equation on $M \setminus S$, then $h$ extends
smoothly as a hermitian metric over the locally free part of $\cE$.
\end{theorem}

Let $A$ be the limiting HYM connection in Theorem~\ref{th:3.1}.
Theorem~\ref{th:3.2} implies that the holomorphic bundle $E|_{M
\setminus S}$ with the holomorphic structure induced by
$\bar{\pt}_{A}$ extends to be a reflexive sheaf $\cE$ over $M$ and
there is a gauge transformation $\s$ on $M \setminus S$ such that
$\s(A)$ extends to a smooth HYM connection on the locally free part of
$\cE$ on $M$. Combining Theorem~\ref{th:3.1} and Theorem~\ref{th:3.2},
the first author was able to describe a compactification of the moduli
space of HYM connections on compact K\"ahler manifolds in
\cite{T}. The idea is to introduce {\em ideal HYM connections} and
define a weak topology on the moduli space of these ideal HYM
connections. Then one can show that the moduli space of the ideal HYM
connections is compact and the moduli space of HYM connections, with
its smooth topology, naturally embeds into it and thus has a natural
analytic compactification. We shall leave the details to the reader as
it is completely analogous to (and simpler than) what we shall do
about vortices and coupled vortices in the next section.
\end{section}

\begin{section}{Compactification of the moduli spaces of vortices and
coupled vortices}
Assume as in Section 2 that $(M, \o)$ is an $m$-dimensional compact
K\"ahler manifold, and that $(E_1, h_1)$ and $(E_2, h_2)$ are Hermitian vector
bundles on $M$. 
\begin{definition}
An {\em ideal coupled ($\t-$) vortex} on the Hermitian bundles $(E_1,
h_1)$ and $(E_2, h_2)$ is a quintuple $(A_1, A_2, \phi, S, C)$ such
that the following holds:

$S$ is a closed subset of finite $H^{n-4}$ measure of $M$. $A_1$ and
$A_2$ are smooth integrable connections on $E_1|_{M \setminus S}$ and
$E_2|_{M \setminus S}$ respectively. $\phi$ is a holomorphic section
of $\Hom(E_2, E_1)|_{M \setminus S}$. $(A_1, A_2, \phi)$ satisfies the
coupled $\tau$-vortex equations (\ref{eq:2.5a}), (\ref{eq:2.5b}) and
(\ref{eq:2.5c}) on $M \setminus S$. $C = (T, \Theta)$ is a $(2m-4)$
dimensional current on $M$. $\supp C = T = \cup_{\a} T_{\a}$ is a
countable union of $(m-2)$ dimensional holomorphic subvarieties of
$M$, and $\frac{1}{8 \pi^2} \Theta |_{T_{\a}}$ are positive
integers. In other words, $\frac{1}{8 \pi^2} C$ is a holomorphic chain
of codimension 2 on $M$. We also require the following energy identity
to be satisfied:
\begin{equation}
\label{eq:4.0}
\YMH_{\t}(A_1, A_2, \phi) +  \| C\| =
\YMH_{\t}(A_1, A_2, \phi) + \int_M \Theta
dH^{2m-4} \lfloor T = E(\tau),
\end{equation}
where $E(\tau)$ is the value of the YMH functional of a smooth coupled
$\tau$-vortex on $(E_1, h_1)$ and $(E_2, h_2)$, which is a geometric
constant depending on $\t$.
\end{definition}
\begin{remark}
Assume that $(A_1, A_2, \phi, S, C)$ is an ideal coupled vortex on
$(E_1, E_2)$. The integrable connections $A_1$ and $A_2$ give
holomorphic structures on $E_1|_{M \setminus S}$ and $E_2|_{M
\setminus S}$. Because of the finiteness of YMH energy of the triple
$(A_1, A_2, \phi)$, we see via Theorem~\ref{th:3.2} that these
holomorphic bundles on $M\setminus S$ extend to reflexive sheaves
$\cE_1$ and $\cE_2$ over $M$. Because $\phi \in \Hom(\cE_1, \cE_2)|_{M
\setminus S}$ and $\Hom(\cE_1, \cE_2)$ is reflexive and hence normal,
we have that $\phi$ extends as a sheaf homomorphism from $\cE_1$ to
$\cE_2$. Therefore, an ideal coupled vortex gives rise to two
reflexive sheaves $\cE_1$ and $\cE_2$ together with a sheaf
homomorphism $\phi \in \Hom(\cE_1, \cE_2)$. It is not hard to see via
Theorem~\ref{th:3.2} that $A_1$ and $A_2$ extend as smooth connections
to the common locally free part of $\cE_1$ and $\cE_2$, which is an
open set whose complement is an analytic subvariety of $M$ of
codimension at least 3.
\end{remark}

We define an equivalence relation among ideal coupled vortices in
the following way. Let $(A_1, A_2, \phi, S, C) \sim (A_1', A_2',
\phi', S', C')$ if $C = C'$ as currents and there exists a closed
subset $S''$ of $M$ with finite $H^{n-4}$ measure such that $(A_1,
A_2, \phi)$ and $(A_1', A_2', \phi')$ are gauge equivalent via smooth
gauge transformations on $M \setminus S''$, i.e., there exist smooth
gauge transformations $g_i$ on $E_i|_{M \setminus S''}$ for $i=1, 2$
such that $(g_1(A_1), g_2(A_2), g_1\circ \phi \circ g_2^{-1}) = (A_1',
A_2', \phi')$ over $M \setminus S''$.

Let the moduli space of ideal coupled $\t$-vortices be given by
$$I\cV_{\t} = I\cV_{\t}(E_1, E_2) = \{ \text{ ideal coupled }
\t\text{- vortices on }(E_1, h_1)\text{ and } (E_2, h_2) \} / \sim. $$
We define the weak topology on $I\cV_{\t}$ via the following notion of
weak convergence. We say a sequence of ideal coupled vortices $(A_{i,
1}, A_{i, 2}, \phi_i, S_i, C_i)$ weakly converges to an ideal coupled
vortex $(A_1, A_2, \phi, S, C)$ if there exist smooth gauge
transformations $g_{i,j}$ on $E_j |_{M \setminus (S_i \cup S)}$ ($j =
1, 2$) such that
$$ (g_{i,1}(A_{i,1}), g_{i,2}(A_{i,2}), g_{i,1}\circ \phi_i \circ
g_{i,2}) \ra (A_1, A_2, \phi) \text{ in } C^{\infty}_c(M \setminus S),
$$
and there is the following convergence of measures,
\begin{equation}
\label{eq:4.0b}
 e_{\t}(A_{i,1}, A_{i,2}, \phi) dv + 8\pi^2 C_i \ra e_{\t}(A_1, A_2, \phi)dv + 8 \pi^2 C. 
\end{equation}
It is easy to see that this weak convergence descends to $I\cV_{\t}$
to define the weak convergence of equivalence classes of ideal coupled
vortices and the limit of a sequence in  $I\cV_{\t}$ is unique. We
note that the moduli space of coupled $\t$-vortices $\cV_{\t}$ with
the smooth topology naturally embeds in $I\cV_{\t}$ if we set the
singularity set $S = \phi$ and the current $C = 0$ for any smooth
coupled vortex  $(A_1, A_2, \phi)$ and thus associates an ideal
coupled vortex to it.

With the above definitions made, we have the following compactness theorem.
\begin{theorem}
\label{th:4.1}
Assume that $(E_1, h_1)$ and $(E_2, h_2)$ are hermitian complex vector
bundles over a compact K\"ahler manifold $M$. Then the moduli space
$I\cV_{\t}$ of ideal coupled vortices on $(E_1, E_2)$ is compact in
the weak topology.
\end{theorem}
\begin{proof}

Let $(A_{i,1}, A_{i,2}, \phi_i, S_i, C_i)$ be a sequence of coupled
$\t$-vortices on hermitian bundles $E_1$ and $E_2$. We adopt the
notations from Section 2. Let $\s$ be defined by $\t$ in
(\ref{eq:2.9}). Let $F = p^*E_1 \oplus p^*E_2 \otimes q^*H^{\otimes
2}$ be the $SU(2)$-equivariant Hermitian vector bundle on $M \times
S^2$.  Each coupled vortex $(A_{i,1}, A_{i,2}, \phi_i)$ corresponds as
in (\ref{eq:2.8a}) to an $SU(2)$-invariant HYM connection
$\tilde{A}_i$ on bundle $F$ with respect to the K\"ahler form
$\o_{\s}$.

We define the concentration set of the sequence $\{ \tilde{A}_i \}$ as
follows (recall definition from (\ref{eq:3.1})),
\begin{equation}
\label{eq:4.1} 
\tilde{S} = \{ x \in M \times S^2 | \lim_{r \ra 0 }  \liminf_{i \ra \infty}
r^{4-(2m+2)} \int_{B_r(x, M \times S^2)} |F_{\tilde{A}_i}|_{\s}^2 dv_{\s} \geq
\e  \},
\end{equation}
where $\e< \e_0$ is to be determined later, $a$ is a constant only
depending on the curvature bound of $M$, and $dv_{\s}$ is volume form
of the metric $\o_{\s}$. Since $\tilde{A}_i$ are $SU(2)$-invariant
connections, the curvature density $ |F_{\tilde{A}_i}|^2$ are
$SU(2)$-invariant functions, hence by the definition $\tilde{S}$ is an
$SU(2)$-invariant subset. Let $\tilde{S} = p^{-1} S$, $S \subset
M$. We call $S$ the concentration set of the sequence $(A_{i,1},
A_{i,2}, \phi_i)$ and we shall show that after gauge transformations,
a subsequence of $(A_{i,1}, A_{i,2}, \phi_i)$ converges to a coupled
vortex $(A_1, A_2, \phi)$ on $M \setminus S$.

We need the following removable singularity theorem. Its proof will be
postponed to the next section.
\begin{theorem}
\label{th:4.1b}
Assume that $(A_1, A_2, \phi, S)$ is an admissible coupled vortex on
trivial vector bundles $E_1$ and $E_2$ over $B_r(0) \subset
\CC^m$. Assume also that the hermitian metrics on the bundles and the
K\"ahler metric on the ball $B_r(0)$ are comparable with the standard
product metrics and K\"ahler metric by a constant factor $c$. Then
there exists $\e_2 = \e_2(c,\t) > 0$, such that if $x \in M$, $0< r
< \injrad(x)$, and
\begin{equation}
\label{eq:4.3c}
r^{4-2m} \int_{B_{r}(x, M)} e(A_{1}, A_{2}, \phi) dv < \e_1,
\end{equation} 
then there exist smooth gauge transformations $g_j$ on $E_j|_{B_{\frac
r2}(x) \setminus S}$ for $j=1, 2$ such that the triple $(g_1(A_1),
g_2(A_2), g_1 \circ \phi \circ g_2^{-1})$ extends smoothly over
$B_{\frac r2}(x)$.
\end{theorem}
We shall also need the following lemma in the proof of Theorem~\ref{th:4.1}.
\begin{lemma}
If $(A_1, A_2,\phi)$ is a coupled $\t$-vortex, then the following
equations hold:
\begin{eqnarray}
\label{eq:4.2a} 
d^*_{A_1} F_{A_1} &=& -2(m-1) J (d_{A_1}(\phi \circ \phi^*)), \\
\label{eq:4.2b} 
d^*_{A_2} F_{A_2} &=& 2(m-1) J( d_{A_2}(\phi^* \circ \phi)),\\
\label{eq:4.2c} 
\bar{\pt}_{A_1 \otimes A_2^*} \phi & =& 0.
\end{eqnarray}
\end{lemma}
\begin{proof}[Proof of Lemma]
Let $\Omega
= \o^{n-2}/(n-2)!$. Then the operator $ * \Omega \wedge :
\La^2_{\CC}(M) \ra \La^2_{\CC}(M)$ has eigenvalues $\pm 1$ and
$\La^2_{\CC}(M) = \La^2(M) \otimes \CC$ decomposes as $\Omega$-self-dual and
$\Omega$-anti-self-dual parts. The space of $\Omega$ self-dual 2-forms
has the decomposition
$$ \La^{+}_{\CC}(M) = \La^0_{\CC}(M) \cdot \omega \oplus
\La^{0,2}_{\CC} \oplus \La^{2,0}_{\CC}. $$
For a connection $A$ on a complex vector bundle $E$ on $M$, we let
$H_A$ be the projection of $F_A$ to the $\La^0_{\CC}(M) \cdot
\omega$ part in the above decomposition. Then we have
$$H_A = (-2i \La F_A) \cdot \omega.$$ Now assume $(A_1, A_2, \phi)$ is
a coupled $\t$-vortex. Since $A_1$ and $A_2$ are
integrable connections, $F_{A_1}$ and $F_{A_2}$ are $(1,1)$ forms. We
have,
\begin{align*} 
d^*_{A_1} F_{A_1} & = - * d_{A_1} * ( F_{A_1} - H_{A_1} +
H_{A_1})\nonumber \\
& = - * d_{A_1} (-\Omega \wedge (F_{A_1} - H_{A_1}) + \Omega \wedge
H_{A_1}) \nonumber \\
& = -2 *(\Omega \wedge d_{A_1} H_{A_1})\quad (\hbox{ by the Bianchi
identity and } d\Omega = 0 ) \nonumber \\
& = -2 *(\Omega \wedge d_{A_1}
(( \phi \circ \phi^* - \t I_{E_1} ) \omega ))\quad (\hbox{ by (\ref{eq:2.5b})} ) \nonumber \\
& = -2 *(\frac{\omega^{m-1}}{(m-2)!} \wedge d_{A_1}(\phi \circ
\phi^*)) \quad (\hbox{ by }d\o = 0)\\
& = -2(m-1)J( d_{A_1}(\phi \circ \phi^*)),
\end{align*}
where $J$ is the complex structure acting on 1-forms. This is exactly
(\ref{eq:4.2a}).  Similarly, we have the equation (\ref{eq:4.2b}) for
$A_2$. (\ref{eq:4.2c}) is just a copy of the vortex equation
(\ref{eq:2.5a}).
\end{proof}

Assume that $x \in M \setminus S$ and $\tilde{x} \in p^{-1}(x) \subset
\tilde{S}$, then by the definition of $S$ and $\tilde{S}$, there
exist $ r \in (0, \dist(x, S)) \cap (0, \injrad (\tilde{x}) )$, and
$N = N(x)>0$, such that for $i \geq N$,
\begin{equation}
\label{eq:4.3} 
r^{4-(2m + 2)}\int_{B_r(\tilde{x}, M\times S^2)}
|F_{\tilde{A}_{i}}|_{\s}^2 dv_{\s} < \e.
\end{equation}
It follows from Lemma~\ref{lem:2.1} that if $r$ satisfies $c(\t) r^4 \leq
\e$ for a suitable constant $c(\t)$, then
\begin{equation}
\label{eq:4.3b}
r^{4-2m} \int_{B_{\frac r2}(x, M)} e(A_{i,1}, A_{i,2}, \phi_i) dv < \frac{2\e}{\s}
\end{equation} 
If we take $\e < 2^{2m-5} \s \e_1$ and fix trivilizations of $E_1$ and
$E_2$ over $B_{r}(x, M)$, then by (\ref{eq:4.3b}) and
Theorem~\ref{th:4.1b}, there exist gauge transformations to make the
triple $(A_{i,1}, A_{i,2}, \phi_i)$ smooth over $B_{\frac r4}(x)$. We
will assume that $(A_{i,1}, A_{i,2}, \phi_i)$ is already in such a
smoothing gauge and hence $\tilde{A}$ is also smooth over
$B_{\frac{r}{4}}(\tilde{x}, M \times S^2)$. Since the left hand sides
of (\ref{eq:4.3}) and (\ref{eq:4.3b}) are scaling invariant, we may
rescale to assume that $r = 8$. In what follows we shall only consider
those $i \geq N$ when we restric our attention on the ball $B_r(x)$.

Since $\tilde{A}$ is now smooth on $B_1(\tilde{x}, \tilde{M})$, the
pointwise a priori estimate for Yang-Mills connections
(Prop.~\ref{prop:3.1}) and (\ref{eq:4.3}) imply that there exists a
uniform constant $C>0$ such that
\begin{equation}
\label{eq:4.4} 
\sup_{B_1(\tilde{x}, M\times S^2)} |F_{\tilde{A}_{i}}|_{\s}^2  \leq C
\e,
\end{equation}
if we assume that $\e \leq \e_0$ for $\e_0$ in Prop.~\ref{prop:3.1}.
 By  Lemma~\ref{lem:2.1} and its proof, we have that $e_{\t}(A_{i,1},
 A_{i,2}, \phi_i) = |F_{\tilde{A}_{i}}|^2 - c(\t) $ and
 $|F_{A_{i,1}}|^2 + |F_{A_{i,1}}|^2 + |d_{A_{i,1} \otimes A_{i,2}^*}
 \phi_i|^2 \leq |F_{\tilde{A}_{i}}|^2$. Thus we have
\begin{eqnarray}
\label{eq:4.5} 
&&\sup_{B_1(x, M)} |F_{A_{i,1}}|^2 + |F_{A_{i,2}}|^2 + |d_{A_{i,1}
\otimes A_{i,2}^*} \phi_i|^2 \leq C \e,\\
\label{eq:4.6} 
&&\sup_{B_1(x, M)} |\phi_i \circ \phi_i^* - \t
I_{E_{1}}|^2 + |\phi_i^* \circ \phi_i + \t' I_{E_2}|^2 \leq C.
\end{eqnarray}
\indent Fix a local trivialization for $E_1$ and $E_2$ on
$B_1(x, M)$. If $\e$ is sufficiently small, by (\ref{eq:4.5}), we may
apply the existence of Coulomb gauges (Theorem 2.7 in Uhlenbeck
\cite{U1}) to find Coulomb gauges for $A_{i,1}$ and $A_{i,2}$ on
$B_1(x, M)$. In other words, there exist gauge transformations
$g_{i,1}$ on $E_1|_{B_1(x)}$ and $g_{i,2}$ on $E_2|_{B_1(x)}$ such
that the connections $A_{i,1}' = g_{i,1}(A_{i,1})$ and $A_{i,2}' =
g_{i,2}(A_{i,2})$ satisify
\begin{eqnarray}
\label{eq:4.7} 
&& d^*A'_{i,j} = 0, \quad \hbox{ on } B_1(x), \quad j=1,2,\\
\label{eq:4.8} 
&& \ast A'_{i,j} = 0, \quad \hbox{ on } \pt B_1(x), \quad j=1,2,\\
\label{eq:4.9}
&& \|A_{i,j}'\|_{L^p(B_1(x))} \leq C'_p \| F_{A_{i,j}}
\|_{L^{\infty}(B_1(x))} \leq C_p \e, \\
&& \qquad \forall 1 \leq p < \infty, j =1,2. \nonumber
\end{eqnarray}
\indent Let $\phi'_i = g_{i,1} \circ \phi \circ g_{i,2}^{-1}$. The equations
(\ref{eq:4.2a}), (\ref{eq:4.2b}), (\ref{eq:4.2c}), (\ref{eq:4.5}), and
(\ref{eq:4.6}) are gauge equivariant and hence they hold if we replace
$(A_{i,1}, A_{i,2}, \phi_i)$ by $(A_{i,1}', A_{i,2}', \phi_i')$. We
observe that (\ref{eq:4.7}) and (\ref{eq:4.8}), (\ref{eq:4.2a}),
(\ref{eq:4.2b}) and (\ref{eq:4.2c}) form an elliptic system for the
triple $(A'_{i,1}, A'_{i,2}, \phi_i')$ over $B_1(x)$. Now
(\ref{eq:4.5}), (\ref{eq:4.6}) and (\ref{eq:4.9}) easily imply that
$\phi_i' \in L^p_1(B_1(x))$ for any $p < \infty$. This together with
(\ref{eq:4.9}) give a starting point to carry out bootstrapping argument and
obtain bounds on the supremum norm of all derivatives of $(A'_{i,1},
A'_{i,2}, \phi_i')$ on $B_1(x)$. In particular, we see that the triple $(A'_{i,1},
A'_{i,2}, \phi_i')$ is smooth on $B_1(x)$ and $(g_{i,1}, g_{i,2})$ are
smooth gauge transformations on $B_1(x)$. The bounds on derivatives of
$(A'_{i,1}, A'_{i,2}, \phi_i')$ are uniform in $i$, hence by passing
to a subsequence, $(A'_{i,1}, A'_{i,2}, \phi_i')$ converges to a
triple $(A'_1, A'_2, \phi')$ in smooth topology on $B_{\frac 12}(x)$.

We may cover the non-concentration set $M \setminus S$ by a countable
union of balls $B_{r_{\a}}(x_{\a})$ such that (\ref{eq:4.3}) and
(\ref{eq:4.3b}) apply with $r = 8r_{\a}$ and $x = x_{\a}$. Applying
the above analysis to each ball $B_{8r_{\a}}(x_{\a})$, by passing to a
subsequence, there exist smooth gauge transformations $(g_{i,1,\a},
g_{i,2,\a})$ on $(E_1, E_2)|_{B_{r_{\a}}(x_{\a})\setminus S_i}$, such
that the sequence $(g_{i,1,\a}(A_{i,1}), g_{i,2,\a}(A_{i,2}),
g_{i,1,\a} \circ \phi_i \circ g_{i,2,\a})^{-1}$ is smooth on
$B_{r_{\a}}(x_{\a})$ and converges in smooth topology to a triple
$(A'_{1,\a}, A'_{2, \a}, \phi'_{\a})$ on $B_{r_{\a}}(x_{\a})$. We can
now use a standard diagonal process of gluing gauges (see for example
4.4.8 in Donaldson and Kronheimer \cite{DK}), again passing to a
subsequence, to obtain smooth gauge transformations $(g_{i,1},
g_{i,2})$ on $M \setminus (S_i \cap S)$, such that $(g_{i,1}(A_{i,1}),
g_{i,2}(A_{i,2}), g_{i,1} \circ \phi_i \circ g_{i,2}^{-1})$ converges
to $(A_1, A_2, \phi)$ in smooth topology on compact subsets of $M
\setminus S$.

We now recall from the dimensional reduction in Section 2 that there
is a one-to-one correspondence between the gauge groups $\cG_1 \times
\cG_2$ and $\cG^{SU(2)}$ and a one-to-one correspondence between
coupled vortices and $SU(2)$-invariant HYM connections. Hence it follows
from the above that after smooth $SU(2)$-invariant gauge transformations,
$\tilde{A_i}$ converges to an $SU(2)$-invariant HYM connection $\tilde{A}$ in
smooth topology on compact subsets of $M \times S^2 \setminus
\tilde{S}$.

Now Prop.~\ref{prop:3.2} on pure Yang-Mills connections
implies that there exists a nonegative density function $\tilde{\Theta}$,
$H^{2m-4}$ measurable on $\tilde{S}$ such that  
\begin{equation}
\label{eq:4.11} 
|F_{\tilde{A}_i}|_{\s}^2 dv_{\s} \ra \mu = |F_{\tilde{A}}|_{\s}^2
 dv_{\s} + \tilde{\Theta} H^{2m-4}\lfloor \tilde{S} \quad \hbox{ as
 measures. }
\end{equation}
We notice that since $|F_{\tilde{A}_i}|_{\s}^2 dv_{\s}$ are
$SU(2)$-invariant measures on $M \times S^2$, the limit measure $\mu$
must also be $SU(2)$-invariant. Recall the definition of the density
function $\tilde{\Theta}$,
\begin{equation}
\label{eq:4.12}
\tilde{\Theta}(\tilde{x}) = \lim_{r \ra 0} \frac{\mu(B_r(\tilde{x}, M
\times S^2))}{r^{2m + 2-4}}, \forall \tilde{x} \in \tilde{S}_b.
\end{equation}
It is clear from the invariance of $\mu$ that $\tilde{\Theta}$ is
$SU(2)$-invariant. Hence $\tilde{\Theta} = \Theta \circ p$
for a function $\Theta$ on $S$. 

By Theorem~\ref{th:3.1} on the blow-up of HYM connections, the blowup
locus of $\tilde{A_i}$ is of the form $\tilde{C} = (\tilde{T},
\tilde{\Theta})$, where $\tilde{T} = \supp \tilde{C}$ is a countable
union of $(m+1-2)$-dimensional holomorphic subvarieties. Because the
current $\tilde{C}$ is $SU(2)$-invariant, we have $\tilde{T} = T
\times S^2 = (\cup_{\a=1}^{\infty} T_{\a}) \times S^2$, where $T_{\a}$
are $(m-2)$-dimensional holomorphic subvarieties of $M$ and the
induced function $\Theta$ satisfies that $(1/ 8\pi^2) \Theta
|_{T_{\a}}$ is a constant positive integer for any $\a$. Let $C' = (T,
\Theta)$, then $\frac{1}{8 \pi^2} C'$ is a holomorphic chain by the
above. Because of the energy identity (\ref{eq:4.0}), the currents
$C_i$ are uniformly bounded in the mass norm, hence passing to a
subsequence, $C_i$ converges to a current $C''$ as
currents. $\frac{1}{8 \pi^2}C''$ is an integral, positive $(m-2, m-2)$
current, hence by the result of King \cite{K}, or Harvey and Shiffman
\cite{HS}, $\frac{1}{8 \pi^2}C''$ is a holomorphic chain. We define $C
= C' + C''$, then $(A_1, A_2, \phi, S, C)$ is an ideal coupled vortex
and the previous argument shows that a subsequence of $(A_{i,1},
A_{i,2}, \phi_i, S_i, C_i)$ weakly converges to $(A_1, A_2, \phi, S,
C)$.
\end{proof}

For hermitian vector bundles $(E_1, h_1)$ and $(E_2, h_2)$ on $M$, as
mentioned before, the moduli space $\cV_{\t}$ of coupled $\t$-vortices
naturally embeds into $I\cV_{\t}$. Hence we have the following
compatification theorem as a corollary of Theorem~\ref{th:4.1}.

\begin{theorem}
\label{th:4.2}
The moduli space $\cV_{\t}$ of coupled $\t$-vortices on hermitian
bundles $(E_1, h_1)$ and $(E_2, h_2)$ over a compact K\"ahler manifold
$(M^m, \o)$ has a compactification $\bar{\cV}_{\t}$ which is embedded
in the space of ideal coupled $\t$-vortices on $M$.
\end{theorem}

In view of the Hitchin-Kobayashi type correspondence established in
Theorem~\ref{th:2.2} and the remark following the definition of ideal
coupled vortices, we have the following corollary.
\begin{theorem}
\label{th:4.3} 
Assume that $(E_1, h_1)$ and $(E_2, h_2)$ are hermitian complex vector
bundles over a compact K\"ahler manifold $(M, \o)$ and $\rank E_2 = 1$
and $\t \notin \cT$ ($\cT$ as defined in (\ref{eq:2.11})). The moduli
space $\gM_{\t}$ of stable holomorphic $\t$-triples on $(E_1, h_1)$
and $(E_2, h_2)$ admits a compactification in the moduli space of
ideal coupled vortices on $(E_1, h_1)$ and $(E_2, h_2)$.
\end{theorem}
\begin{remark}
From the remark following the definition of ideal coupled vortices (at
the beginning of this section), we know that an ideal coupled vortex
$(A_1, A_2, \phi, S, C)$ gives rise to reflexive sheaves $\cE_1$ and
$\cE_2$ and a sheaf homomorphism $\phi \in \Hom(\cE_1, \cE_2)$. We
should be able to define suitable notions of stability and
semi-stability for such triples $(\cE_1, \cE_2, \phi)$ (by algebraic
criteria) so that the stability of a triple is equivalent to the
existence of a singular coupled vortex $(A_1, A_2, \phi)$ on it. This
will be the Hitchin-Kobayashi correspondence for stable triples in the
sheaf version. With this work done, Theorem~\ref{th:4.3} could be
written in a nicer way, i.e., the compactification of the moduli space
of stable bundle triples lies in the moduli space of stable sheaf
triples, which will be in the algebraic category.
\end{remark}

Now assume that $(E, h)$ is a hermitian vector bundle over $M$. In order to
get a compactification of the moduli space $V_{\t}(E)$ of
$\t$-vortices on $(E, h)$, we use the embedding of $V_{\t}(E)$ into
$\cV_{\t}(E, L)$, the moduli space of coupled $\t$-vortices on $E$ and
$L$, where $L$ is the trivial smooth line bundle over $M$
with the product metric $h_0$. We have the following compactification
theorem of the moduli space of vortices.
\begin{theorem}
\label{th:4.4}
The moduli space $V_{\t}$ of $\t$-vortices on a hermitian vector
bundle $E$ over a compact K\"ahler manifold $(M, \o)$ admits a
compactification in the space of ideal coupled $\t$-vortices on $E$
and $L$, where $L$ is the trivial line bundle with the product metric
on $M$.
\end{theorem}
\begin{remarks}
1) Assume that $(A_1, A_2, \phi, S, C)$ is an ideal coupled
$\t$-vortex on $E$ and $L$. Then it gives rise to reflexive sheaves
$\cE_1$ and $\cL$.  $\cL$ is a line bundle because it is a rank 1
reflexive sheaf. The degree of $\cL$ is again equal to $0$, that of
the trivial bundle, because in passing to a weak limit, the
topological change happens only on the blowup set, a codimension 2
set, and that doesn't affect the first chern class; hence $\cL$ is
topologically trivial. Now we may determine the connection $A_2$
uniquely in terms of the connection $A_1$, $\phi$ and the holomorphic
structure $\cL$ from the equation (\ref{eq:2.5c}). Thus the ideal
coupled vortex essentially can be viewed as an `ideal vortex' $(A_1,
\phi, S, C)$ on the bundle $\cE_1 \otimes \cL^*$. And the boundary
points of the compactification $\bar{V_{\t}}$ are actually `ideal
vortices'. However, we would not use the concept of `ideal vortices'
here as it would not make things much simpler. \\ 2) The
compactification given in Theorem~\ref{th:4.4} is somewhat indirect
since in order to compactify vortices, we turn to ideal {\em coupled}
vortices. However, we do not know a more direct way of compactifying
the space of vortices. Notice in particular the blow-up set of a
sequence of vortices $(A_i, \phi_i)$ is not where the YMH energy of
$(A_i, \phi_i)$ concentrates. Instead it is the set where the YMH
energy of the sequence of coupled vortices $(A_i, A_i', \phi_i)$
concentrates, where $A_i'$ is a connection on $L$ determined by $(A_i,
\phi_i)$ from (\ref{eq:2.5c}). In this regard, we may say that the
blow-up phenomena of vortices is only clear when we put them in the
setting of coupled vortices.
\end{remarks}
\indent Via the Hitchin-Kobayashi type correspondence given by
Theorem~\ref{th:2.1}, we have,
\begin{theorem}
\label{th:4.5} 
Assume that $\t \notin \cT$ ($\cT$ as defined in
(\ref{eq:2.4b})). Then the moduli space $M_{\t}$ of stable $\t$-pairs
on a Hermitian complex vector bundle $(E, h)$ over compact a K\"ahler
manifold $(M^m, \o)$ admits a compactification in the moduli space of
ideal coupled vortices on $E$ and $L$, where $L$ is the trivial bundle
with the product metric on $M$. 
\end{theorem}
\begin{remarks}
1) Again, like in the remark following Theorem~\ref{eq:4.4}, we should
be able to define an algebraic concept of stable (or semi-stable)
pairs $(\cE, \phi)$ of a coherent reflexive sheaf and a section of
it, and to interpretate the compactification in the algebraic
category. We hope to clarify these issues in a future paper. \\
2) When $M$ is a K\"ahler surface, let $(A_1, A_2, \phi, C)$ be an
   ideal coupled vortex on hermitian bundles $(E_1, h_1)$ and $(E_2,
   h_2)$. Then because reflexive sheaves are locally free on K\"ahler
   surfaces, we find that $(A_1, A_2, \phi)$ extends as a smooth
   coupled vortex on some hermitian bundles $(E_1', h_1')$ and $(E_2',
   h_2')$ on $M$. $C$ is now given by a finite set of points with
   multiplicities on $M$. Hence the compactification is much easier to
   describe and it is reminiscent of the well-known moduli spaces of
   self-dual connections on four-manifolds. It would be a natural
   question to ask whether we can describe the compactified moduli
   spaces more specificly, and in particular, compute their topological
   invariants. \\
3) In order to understand the compacitification
   better, we have the following question. Are there any effective
   bounds on the measure of the singular set of a reflexive sheaf in terms of
   topological or analytical quantities?  \\
\end{remarks}
\end{section}

\begin{section}{Removable singularity theorems}
In this section we give the proofs of Theorem~\ref{th:1.4},
\ref{th:1.5} and \ref{th:4.1b}. We shall first prove the following
$\e$-regularity theorem for HYM connections.
\begin{theorem}
\label{th:5.1}
Assume that $E$ is a trivial complex vector bundle over $B_{r}(0)
\subset \CC^m$ with the product hermitian metric. Assume also that the
hermitian metric on the bundle and the K\"ahler metric on the ball
$B_r(0)$ are comparable with the standard product metrics and K\"ahler
metric by a constant factor $c$. Then there exists a consant $\e_1(m, c) > 0$,
such that if $(A, S)$ is an admissible HYM connection on $E$, with
\begin{equation}
\label{eq:5.1}
r^{4-2m}\int_{B_r(0)} |F_A|^2 dv \leq \e_1, 
\end{equation}
then there exists a smooth gauge transformation $\s$ on $B_{\frac
r2}(0) \setminus S$, such that $\s(A)$ can be extended smoothly over $B_{\frac
r2}(0)$.
\end{theorem}
\begin{proof}
Because (\ref{eq:5.1}) is scaling invariant, we may rescale and assume that
$r = 1$. Let $E_0$ be the holomorphic bundle over $B_r(0) \setminus S$
which is topologically $E$ over $B_r(0) \setminus S$ and has the
holomorphic structure given by $\bar{\pt}_A$. We first resort to
Theorem~\ref{th:3.2} to extend $E_0$ to a reflexive sheaf $\cE$ over
$B_r(0)$. Let $S_1$ be the singular set of the reflexive sheaf $\cE$,
i.e. the subset of $M$ where $\cE$ is not locally free. It is a
standard fact that $S_1$ is an analytic subvariety of dimension at
most $m-3$.  Theorem~\ref{th:3.2} also implies that the hermitian
metric $h$ extends smoothly over the locally free part $B_r(0)
\setminus S_1$. Hence the connection $A$, determined by $h$ and the
holomorphic structure, also extends smoothly over $B_r(0) \setminus
S_1$.

Since $S_1$ is stratified by smooth complex submanifolds of $M$, we
shall make induction on the complex dimension of the top strata $S_0$
of $S_1$ to show that, under the assumption (\ref{eq:5.1}), $S_1 =
\emptyset$. 

If $S_0 = \emptyset$, then $S_1 = \emptyset$. Assume that we have
shown that if $\dim S_0 < k$, then $S_1 = \emptyset$. Assume now that $\dim
S_0 = k$. Take a generic point $x_0 \in S_0$. With a suitable choice of
coordinates, a neighborhood of $x_0$ in $B_1(0)$ can be written in the
form of $N = B^{m-k}_s \times B^k_s$ such that $x_0 = (0, 0)$, and
$S_1 \cap N = S_0 \cap N = \{0\} \times B^k_s$, where $B^k_s$ stands
for the ball of radius $s$ centered at origin in $\CC^k$. Again after
rescaling, we may assume that $s = 2$.

We first claim that the vector bundle $\cE|_{N \setminus S_0} =
\cE|_{(B^{m-k}_1 - \{0\}) \times B^k_1 }$ is trivial as a smooth
vector bundle over $N \setminus S_0$. Let $F = \cE|_{(B^{m-k}_{\frac
12}-\{0\}) \times \{0\}}$ and $A' = A|_{\pt B^{m-k}_{\frac 12} \times
\{0\}}$.  If $\e_1 < \e_0$, we have by the a priori estimates
(Prop.~\ref{prop:3.2}) that
$$ |F_{A'}|^2(y) \leq |F_A|^2(y) \leq C(m) \e_1, \quad \forall y \in
\pt B^{m-k}_{\frac 12} \times \{0\} $$
We can then apply the argument of Lemma 2.2 of Uhlenbeck \cite{U1} to
assert that the bundle $F$ restricts to a trivial smooth bundle over
$\pt B^{m-k}_1 \times \{0\}$, if $\e_1$ is sufficiently small. Let $\{
e_1, \ldots, e_l \}$ be a frame of $F$ over $\pt B^{m-k}_1 \times
\{0\}$ and fix a smooth connection on $\cE|_{N \setminus S_0}$, we
may use parallel transport to obtain a smooth trivialization $\{e_1,
\ldots, e_n \}$ of $\cE|_{N \setminus S_0}$ and the claim is proved.

The trivialization of $\cE$ away from the singular set allows us to
express $A$ as a matrix valued $1$-form. If $\e_1$ is sufficiently
small, the singularity removal theorem of Tao and Tian \cite{TT}
implies that there exists a smooth gauge transformation $g$ on
$(B^{m-k}_{\frac 12} -\{0\}) \times B^{k}_{\frac 12} $ such that
$g(A)$ extends as a smooth HYM connection over a trivial bundle $H$
over $B^{m-k}_{\frac 12} \times B^k_{\frac 12}$. Define the
holomorphic structure on $H$ by $\bar{\pt}_{g(A)}$. Let $\{ f_1,
\ldots, f_l \}$ be a frame of holomorphic sections of $H$ over
$B^{m-k}_{\frac 12} \times B^{k}_{\frac 12} $. Since $\bar{\pt}_{g(A)}
= g \circ \bar{\pt}_A \circ g^{-1}$, $\{ g^{-1} f_1, \ldots, g^{-1}
f_l \}$ gives a holomorphic frame for the reflexive sheaf $\cE$ over
$(B^{m-k}_{\frac 12} -\{0\}) \times B^{k}_{\frac 12} $. Since $\cE$ is
a reflexive coherent sheaf and $\cE$ agrees with the trivial
holomorphic bundle away from a singular set of codimension at least
$3$, it follows that $\cE$ is actually a trivial holomorphic bundle
over $B^{m-k}_{\frac 12} \times B^k_{\frac 12}$. This implies that $\cE$ is
locally free at $x_0$. It is a contradiction and the claim is
established.
\end{proof}

\begin{proof}[Proof of Theorem~\ref{th:4.2}]
We shall use the notations from Section 2. Let $\tilde{A}$ be the
$SU(2)$-invariant admissible HYM connection on the $SU(2)$-invariant
bundle $F = p^*E_1 \oplus (p^*E_2 \otimes q^* H^{2})$ determined by
$(A_1, A_2, \phi)$ as in Section 2. If $\e_2$ is sufficiently small,
we apply Theorem~\ref{th:5.1} to $\tilde{A}$ over local patches of
$B_r(0) \times S^2$ and see that locally $\tilde{A}$ can be
extended to a smooth HYM connection after a gauge transformation and
the holomorphic bundle $F|_{(B_r(0) \setminus S)\times S^2}$
locally extends across the singularity to a holomorphic bundle. Since
we know from Bando and Siu's theorem (Theorem~\ref{th:3.2}) that
$F|_{(B_r(0) \setminus S)\times S^2}$ extends uniquely as a
reflexive sheaf $\cF$ over $B_r(0) \times S^2$, it follows that
$\cF$ is a holomorphic bundle over $B_{r}(0) \times S^2$ with a
smooth Hermitian-Einstein metric $h$. $\cF$ and its metric are
$SU(2)$-invariant because so are $F$ and its metric.  Let $\cE_1$,
$\cE_2$ be the reflexive sheaves extending $E_1$, and $E_2$ over
$B_r(0)$ (by Theorem~\ref{th:3.2}). Then $\cF$ has a splitting
$$ \cF = p^*\cE_1 \oplus (p^*\cE_2 \otimes q^*H^2),$$
because the splitting is valid on $(B_{\frac r2}(0) \setminus S)
\times S^2$ and the sheaves involved are reflexive. With the following lemma,
we deduce that $\cE_1$ and $\cE_2$ are smooth holomorphic
bundles over $B_{\frac r2}(0)$. Since the hermitian metric $h$ on $\cF$
is $SU(2)$-invariant, it
must be of the form
$$ h = \begin{pmatrix}
                \tilde{h}_1 & 0 \\
                0  & \tilde{h}_2
        \end{pmatrix},
$$
where $\tilde{h}_1 = p^*\bar{h}_1 $, $\tilde{h}_2 = p^*\bar{h}_2
\otimes q^* h_2'$, and $\bar{h}_i$ is a metric on $\cE_i$ and $h_2'$
is the $SU(2)$-invariant metric on $H^2$ (see Prop. 3.2 of
\cite{G1}). $\bar{h}_i$ and the holomorphic structures on $\cE_i$ then
give a smooth coupled $\t$-vortex $(A_1', A_2', \phi')$ on $B_{\frac
r2}(x, M)$. $(\cE_i, \bar{h}_i)|_{B_{\frac r2}(0)}$ is
isomorphic to $(E_i, h_i)|_{B_{\frac r2}(0)}$ as hermitian holomorphic
vector bundles because by definition, the former ones are extensions of
the later ones. It follows that $(A_1', A_2', \phi)$ is  gauge
equivalent to $(A_1, A_2, \phi)$ on $B_{\frac r2}(0)$ and gives the
desired smooth extension.
\end{proof}
\begin{lemma}
Assume that $\cE_1$ and $\cE_2$ are two coherent sheaves on a complex
manifold $M$ such that $\cE_1 \oplus \cE_2 = \cF$, where $\cF$ is a
holomorphic vector bundle of finite rank over $M$. Then $\cE_1$ and $\cE_2$ are locally free.
\end{lemma}
\begin{proof}[Proof of lemma]
The statement is local in nature. Let $x$ be any point of $M$. Because of
$(\cE_1)_x \oplus (\cE_2)_x = \cF_x$, we see that $(\cE_i)_x$ ($i = 1,
2$) are projective, hence are free. Because the sheaves are coherent,
we see that $\cE_1$ and $\cE_2$ are free in a neighborhood of
$x$. That finishes the proof of the lemma.  
\end{proof}
In the end, we remark that the proofs of Theorem~\ref{th:1.4} and
\ref{th:1.5} follow easily from Theorem~\ref{th:5.1} and \ref{th:4.1b}
through a standard covering argument and we shall leave the proof to
the reader.
\end{section}

\vspace{5ex}
{\centering Department of Mathematics, Massachussetts Institute of Technology, 77 Mass. Ave.,
Cambridge, MA 02139, USA\\
email: tian@math.mit.edu\\
\vspace{2ex}
Department of Mathematics, Stanford University, Stanford, CA
94305-2125, USA\\
email: byang@math.stanford.edu\\}

\end{document}